\documentclass[english]{amsproc}
\usepackage{amsmath,amssymb,amsthm,amsfonts,amscd,mathrsfs}

\usepackage{babel}
\usepackage[all]{xy}
\usepackage{tikz}
\usetikzlibrary{decorations.pathmorphing}
\usepackage{tikz-cd}
 \usepackage{graphicx} 
 \usepackage{epstopdf}
 \usepackage{pb-diagram}
 \usepackage[]{todonotes}
\usepackage{mathtools}

\newcommand{\bdm}{\begin{displaymath}}
\newcommand{\edm}{\end{displaymath}}

\newcommand{\Z}{\mathbb{Z}}

\newcommand{\co}{\colon\thinspace}
\newcommand{\TC}{\mathbf{TC}}

\newcommand{\TCG}{\mathbf{\TC}_G}

\newcommand{\nil}{\mbox{\rm nil}}
\newcommand{\id}{\operatorname{id}}
\newcommand{\ev}{\operatorname{{\it ev}}}
\newcommand{\cat}{\operatorname{cat}}

\newcommand{\secat}{\operatorname{secat}}

\newtheorem{defn}{Definition}[section]

\newtheorem{exam}[defn]{Example}

\theoremstyle{plain}
\newtheorem{thm}[defn]{Theorem}
\newtheorem{prop}[defn]{Proposition}
\newtheorem{lemma}[defn]{Lemma}

\begin{document}

\title[Equivariant topological complexities]{Equivariant topological complexities}

\author{Andr\'es \'Angel}
\address{Departamento de Matem\'aticas, Universidad de los Andes, Carrera 1 N. 18A - 10 Bogot\'a, Colombia}

\email{ja.angel908@uniandes.edu.co}

\author{Hellen Colman}
\address{Department of Mathematics, Wright College,
 4300 N.\ Narragansett Avenue, Chicago, IL 60634 USA}
\email{hcolman@ccc.edu}
\thanks{The work of H.C. was supported in part by Simons Foundation. The work of A.A. was supported in part by the FAPA funds from Vicerrector\'ia de Investigaciones de la Universidad de los Andes}

\date{\today}

\begin{abstract}
The aim of this article is to review different generalizations of the the notion of topological complexity to the equivariant setting. In particular, we review the relation (or non-relation)  between these notions and the topological complexity of the quotient space and the topological complexity of the fixed point sets. We give examples of calculations and stress the question: When the action is free, do we recover the topological complexity of the quotient?  
\end{abstract}

\maketitle

\section{Introduction}
The topological complexity of spaces with
group actions was first introduced by Colman and Grant \cite{Colman2012} as a navigational complexity quantifier of certain mechanical problems best described by groups acting on spaces and on the other hand as a tool for obtaining
useful information in classical nonequivariant topological complexity. 

Since then there have been several approaches to defining other equivariant versions of topological complexity emphasizing in different degrees one or the other of these two broad objectives. Namely,

\begin{itemize}

\item (Colman-Grant) Equivariant topological complexity.
\item (Lubawski-Marzantowicz) Invariant topological complexity.
\item (Dranishnikov) Strongly equivariant topological complexity.
\item (B{\l}aszczyk-Kaluba) Effective topological complexity.
\end{itemize}

In this paper we give a general  view of the different definitions as well as a summary of their properties.

Most results considered in our exposition are based on the original articles by the respective authors \cite{Colman2012, Lubawski2014, Dranishnikov2015,effective} and calculations from \cite{BlaszczykKaluba,productprojective}. We also include some new results.

The layout of the paper is as follows. In section 2 we review some general equivariant notions and we define the equivariant Lusternik-Schnirelmann category of a group action. In section 3 we introduce the Farber's topological complexity of a topological space and discuss some of its properties.  Section 4 is the core of the article and provides a survey  of the four definitions of equivariant versions of topological complexity and its properties. In section 5 we finish with comments on the role of the notion of Morita equivalence of group actions in the study of topological complexity of spaces with symmetries. We propose the study of Morita invariance for the {\em invariant topological complexity}. Moreover, we aim to develop a theory of topological complexity for  orbifolds \cite{AngelColman} which will generalize and give as particular cases the classical topological complexity of spaces and a topological complexity for group actions.  In general, an {\em orbifold topological complexity} will give a version of topological complexity for all contexts modeled by classes of equivalence of group actions such as topological spaces, group actions, foliations and stacks.

\section{Equivariant Notions}

In this paper, $G$ will always denote a compact Hausdorff topological group acting continuously on a Hausdorff space $X$ on the left. In this case, we say that $X$ is a {\em $G$-space}. For each $x\in X$ the {\em isotropy group} $G_x=\{h\in G\mid  hx=x\}$ is a closed subgroup of $G$. The set $Gx=\{gx\mid g\in G\}\subseteq X$ is called the {\em orbit} of $x$. 

The {\em orbit space} $X/G$ is the set of equivalence classes determined by the action, endowed with the quotient topology. Since $G$ is compact and $X$ is Hausdorff, $X/G$ is also Hausdorff, and the {\em orbit map} $p\co X\to X/G$ sending a point to its orbit is both open and closed.

 If $H$ is a closed subgroup of $G$, then $X^H=\{x\in X|\; hx=x \mbox{ for all }h\in H\}$ is called the {\em $H$-fixed point set} of $X$. We call $x$ a {\em global fixed point} if $x \in X^G$.

A $G$-space $X$ is said to be {\em $G$-connected} if
the $H$-fixed point set $X^H$ is path-connected for every closed subgroup $H$ of
$G$.

\subsection{G-homotopy}

Let $X$ and $Y$ be $G$-spaces.
Two $G$-maps $\phi, \psi\co  X\to Y$ are {\em $G$-homotopic}, written $\phi\simeq_G \psi$, if there is a $G$-map $F\co X \times I \to Y$ with $F_0=\phi$ and $F_1=\psi$, where $G$ acts trivially on $I$ and diagonally on $X\times I$.

If there exist $G$-maps $\phi\co X\to Y$ and $\psi\co Y\to X$ such that $\phi\psi\simeq_G \id_Y$ and $\psi\phi\simeq_G \id_X$, then $\phi$ and $\psi$ are {\em $G$-homotopy equivalences}, and $X$ and $Y$ are {\em $G$-homotopy equivalent}, written $X\simeq_G Y$.

 \begin{defn} We say that a $G$-invariant subset $U \subseteq X$ is {\em $G$-compressible into a $G$-invariant subset $A \subseteq X$}, 
if the inclusion map $i_U : U \rightarrow X$ is $G$-homotopic to a $G$-map $c: U \rightarrow X$ with $c(U) \subseteq A$.
\end{defn}

 \begin{defn} A $G$-invariant subset $U \subseteq X$  is called {\em $G$-categorical} if $U$ is $G$-compressible into a single orbit. 
\end{defn}
We say that the space $X$ is {\em $G$-contractible} if $X$ is $G$-categorical.

The following two lemmas emphasize the importance of 
 $G$-connectedness and the existence of global fixed points in relation to $G$-homotopy. They 
 will be repeatedly used in many subsequent arguments.
 
\begin{lemma}[Conservation of isotropy]\
Let $X$ be a $G$-connected $G$-space, and let $x, y \in X$ such that $G_x\subseteq G_y$. Then there exists a $G$-homotopy $F\co Gx \times I \to X$ such that $F_0=i_{Gx} : Gx \hookrightarrow X$ and $F_1(Gx)\subseteq Gy$.
\end{lemma}

\begin{lemma}[$G$-categorical is $G$-compressible into a point]\
Let $X$ be a $G$-connected $G$-space with $X^G \neq \emptyset$ , and let $x \in X^G$. Then every $G$-categorical subset $U$ of $X$ is $G$-compressible to $x$. \label{Gcat}
\end{lemma}

\subsection{Equivariant Lusternik-Schnirelmann category}

Marzantowicz \cite{Marzantowicz} studied the equivariant LS-category as a generalization of the usual notion of Lusternik-Schnirelmann category adapted to the equivariant context.

 \begin{defn}
Given a $G$-space $X$, the {\em equivariant LS-category}, $\cat_G(X)$, of $X$ is the least integer $k$ such that $X$ may be covered by $k$ $G$-invariant open sets   $\{U_1,\ldots ,U_k\}$, each one $G$-categorical.
\end{defn}

By definition we have that $\cat_G(X)=1$ iff $X$ is $G$-contractible. Like in the nonequivariant version, if $M$ is a manifold, we have that $\cat_G(M)$ is a lower bound for the number of critical points of invariant differentiable functions $f:  M \rightarrow \mathbb{R}$.

\begin{prop} Let $X$ and $Y$ be  $G$-spaces.
\begin{enumerate}
\item If $X\simeq_G Y$ then $\cat_G(X)=\cat_G(Y)$.
\item If $\cat_G(X)=1$ and $X^G \neq \emptyset$ then $X$ is $G$-compressible to a point.
\item $|\pi_0(X^G)| \leq \cat(X^G) \leq \cat_G(X)$.
\item $\cat(X/G) \leq \cat_G(X)$.
\end{enumerate}
\end{prop}

If $G$ is a compact Lie group, Marzantowicz proved that the equivariant LS-category of a free space is just the LS-category of the quotient. The main ingredient of the proof is the use of the Palais covering homotopy theorem \cite{Palais} which gives conditions to find equivariant lifts of homotopies between quotient spaces.
\begin{thm}\cite{Marzantowicz} If  a compact Lie group $G$ acts freely on a metrizable space  $X$, then 
$$
\cat_G(X) =\cat(X/G).
$$
\end{thm} \label{Marza}
More generally, the result holds if  $X$ has one orbit type. 

For example, from this theorem follows the calculation of the equivariant LS-category of the antipodal action 
of $\mathbb{Z}_2$ on $S^n$:
$$
\cat_{\mathbb{Z}_2}(S^n) = \cat(\mathbb{RP}^n)=n+1.
$$

A cohomological lower bound for $\cat_G(X)$ is given in \cite{Marzantowicz} using the nilpotency of reduced equivariant cohomologies. An equivariant cohomology theory is called {\em singular } if for every closed subgroup $H \subseteq G$ we have that  $H^*_G(G/H)=0$ for $*>0$.

\begin{thm}\cite{Marzantowicz} Let $H^*_G(\cdot)$ be a singular multiplicative $G$-cohomology theory. If we denote by $\mathfrak{Z}$  the reduced equivariant cohomology $\tilde{H}^*_G(X)$, then $\cat_G(X) 
\geq  \nil{(\mathfrak{Z})}$.
\end{thm}

For any group $G$ we have the universal example of free $G$-space $EG$. It is a free $G$-space that is contractible and such that the projection to the quotient space $BG=EG/G$ is a principal $G$-bundle. When $G$ is a finite group we have that $\cat(BG)=\infty$ because $BG$ has infinite cohomological dimension.

When  $G$ is an infinite discrete group, Palais covering homotopy theorem cannot be applied. Instead, by using the homotopy lifting property of covering spaces we have a similar result as in the previous theorem:
$$\cat_G(EG)=\cat(EG/G)=\cat(BG)$$
which is equal to $cd(G)+1$ where $cd(G)$ is the cohomological dimension of $G$.

A product inequality $\cat_G(X\times Y) \leq \cat_G(X)+\cat_G(Y)-1$ is treated in theorem 3.15 in \cite{Colman2012}, but the theorem as stated is not true. A counterexample and the correct hypotesis for the result to be true, are discused in theorem 2.23 and example 6.4 in \cite{Bayeh}.

The assumptions that $X$ is $G$-connected and $X^G\neq\emptyset$ are enough to have that  $\cat_G(X \times X) \leq 2\cat_G(X)-1$.

We have also the generalization to the equivariant case of the Clapp-Puppe category \cite{Clapp} that is relevant to the study of invariant topological complexity and strongly equivariant topological complexity. This generalization is defined in \cite{Marzantowicz}.

\begin{defn}
Let  $A\subseteq X$ be a $G$-invariant subset of $X$. The {\em equivariant $A$-category},  $_{A}cat_G(X)$,  is the least integer $k$ such that $X$ may be covered by $k$ $G$-invariant open sets   $\{U_1,\ldots ,U_k\}$, each  $G$-compressible into $A$.
\end{defn}

If $* \in X^G$ then $cat_G(X) \leq _{\{*\}}cat_G(X)$ and if $X$ is $G$-connected, 
using lemma \ref{Gcat}, we have that every $G$-categorical set is $G$-compressible to the point $*$ and thus
$_{\{*\}}cat_G(X) \leq cat_G(X)$. Note that in this case we have 
$$
\cat(X) \leq {_{\{*\}}cat_G(X)} = cat_G(X).
$$


\section{Topological complexity}
Let $X$ be the space of configurations of a mechanical system. A {\em motion planning algorithm} is a set of rules that to each pair of configurations (initial and final) assigns a path  between them.

Let $X^I$ be the space of all paths in a space $X$. Consider the evaluation map $\ev \colon X^I \to X \times X$ given by
\[ \ev(\gamma) = \big(\gamma(0), \gamma(1) \big).\]

A \textit{motion planner} on an open subset $U \subseteq X \times X$ is a section of the evaluation map $\ev$ over $U$, i.e. a (continuous) map $s: {U}\to X^I$
  such that  the following diagram commutes:

$$\xymatrix{& X^I \ar[d]^{\ev} \\ {U}  \ar[ur]^{s} \ar[r]^{}& X\times X }$$

The {\em topological complexity} of a space $X$, denoted $\TC(X)$, is the least integer $k$ such that there exists an open cover of $X\times X$ by $k$ open sets on each of which there is a motion planner \cite{Farber2003a, Farber2004}.

\begin{prop} Let $X$ and $Y$ be spaces.
\begin{enumerate}
\item If $X\simeq Y$ then ${\TC}(X)={\TC}(Y)$.
\item ${\TC}(X)=1$ iff $X$ is contractible.
\item ${\TC}(X)=\infty$ if $X$ is not connected.
\end{enumerate}
\end{prop}

\begin{thm} Let $X$ be a connected space. 
\begin{enumerate}
\item $\cat(X) \leq {\TC}(X) \leq \cat(X \times X)$.
\item ${\TC}(X) \leq 2 \dim X+1$ where $\dim$ denotes the covering dimension.
\end{enumerate}
\end{thm}

Topological complexity can be defined equivalently in terms of the sectional category of a map \cite{Schwarz} and in terms of the Clapp-Puppe category:

\begin{thm} For a space $X$,   the following statements are equivalent: 
\begin{enumerate}
\item ${\TC}(X) \leq n$.
\item $\secat(ev) \leq n$: there exist open sets $U_1,\ldots, U_k$ which cover $X \times X$ and sections $s_i :  U_i \rightarrow X^I$ such that $ev \circ s_i$ is homotopic to the inclusion map $U_ i \hookrightarrow X \times X$.
\item $_{\Delta(X)}\cat(X \times X) \leq n$: there exist open sets $U_1,\ldots, U_k$ which cover $X \times X$ such that each is compressible into $\Delta(X)$.

\end{enumerate}
\end{thm}

Cohomological lower bounds are given by the nilpotency of  the kernel of zero divisors of $X$. Let $H^*(-)$ be the cohomology with coefficients in a field and $\mathfrak{Z}=\ker(\cup) $ where
$
\cup: H^*(X) \otimes H^*(X) \rightarrow H^*(X)
$ is the cup product homomorphism.

\begin{prop}

${\TC}(X) > \nil{(\mathfrak{Z})}$ where $\mathfrak{Z}$ is the kernel of zero divisors.

\end{prop}

\begin{exam}{Topological complexity of spheres is given by:}
\begin{enumerate}
\item ${\TC}(S^0)=\infty$.
\item ${\TC}(S^n)=2$ for $n$ odd.
\item ${\TC}(S^n)=3$ for $n> 0$ even.
\end{enumerate}
\end{exam}

\begin{thm}
For the trivial $F$-bundle $B \times F \rightarrow B$, we have that
$$
{\TC}(B \times F) \leq {\TC}(B) + {\TC}(F) -1.
$$
\end{thm}


\section{Equivariant versions of TC}
We review in this section the different approaches to define a topological complexity in the context of group actions. From now on, $X$ will be always a $G$-space. All four versions reduce to the ordinary (nonequivariant) $\TC(X)$ when the action of $G$ on $X$ is trivial. Moreover, all four versions are invariant under $G$-homotopy type.

\subsection{Equivariant topological complexity (Colman-Grant)}
The evaluation map $\ev\co X^I\to X\times X$ is a $G$-fibration with respect to the actions
$$G\times X^I\to X^I, \qquad G\times (X\times X)\to X\times X,$$
$$g(\gamma)(t)=g(\gamma(t)),\qquad g(x,y)=(gx,gy).$$
\begin{defn}
The {\em equivariant topological complexity} of $X$, $\TCG(X)$, is the least integer $k$ such that
  $X \times X$ may be covered by $k$ $G$-invariant open sets
  $\{U_1,\ldots ,U_k\}$, on each of which there is a $G$-equivariant map 
   $ s_i:{U}_i\to X^I$ such that  the  diagram commutes:
$$\xymatrix{& X^I \ar[d]^{\ev} \\ {U}_i  \ar[ur]^{s_i} \ar[r]^{}& X\times X }$$
\end{defn}
In other words, the equivariant topological complexity of a $G$-space $X$  is 
the minimum number of $G$-invariant open sets needed to cover $X\times X$, on each of which  the free path fibration $\ev$ admits a local $G$-equivariant section. If no such integer exists then we set $\TCG(X)=\infty$.


\begin{prop} Let $X$ and $Y$ be  $G$-spaces. 
\begin{enumerate}
\item For a $G$-connected space $X$  with $X^G \neq \emptyset$, $\TC_G(X)=1$ iff $X$ is $G$-contractible.
\item $\TC_G(X)=\infty$ if $X$ is not $G$-connected.
\end{enumerate}
\end{prop}

 Moreover there are inequalities relating $\TCG(X)$ to the equivariant and non-equivariant topological complexities of the various fixed point sets.

\begin{prop} \label{prop3} If $H$ is a closed subgroup of $G$, then
\begin{enumerate}
\item ${\TC}(X^H) \leq \TC_G(X)$.
\item $\TC_H(X) \leq \TC_G(X)$, in particular 
${\TC}(X) \leq \TC_G(X)$.
\end{enumerate}
\end{prop}

The next results describe the basic relationship of equivariant topological complexity with equivariant Lusternik-Schnirelmann category.

\begin{prop} Let $X$ be a $G$-connected space.

\begin{enumerate}
\item $\TC_G(X) \leq \cat_G(X \times X)$.
\item If $X^G\neq\emptyset$, then 
$\cat_G(X) \leq  \TC_G(X) $.
\end{enumerate}
\end{prop}

The assumptions that $X$ is $G$-connected and $X^G\neq\emptyset$ are enough to have that  $\cat_G(X \times X) \leq 2\cat_G(X)-1$ and it follows that $\TC_G(X) \leq 2\cat_G(X)-1$.

Since an invariant open subset $U\subseteq
X\times X$ has a $G$-equivariant section 
   $ {U}\overset{s}{\to}X^I$ of
$\ev\co X^I\to X\times X$ if and only if the inclusion $i_U\co U\hookrightarrow X\times X$
is $G$-homotopic to a map with values in the diagonal
$\Delta(X) \subseteq X\times X$, we can reformulate the definition of equivariant topological complexity in terms of equivariant deformations to the diagonal. See lemma 3.5 in \cite{Lubawski2014} for the following:

\begin{thm} For a $G$-space $X$, the following statements are equivalent: 
\begin{enumerate}
\item ${\TC}_G(X) \leq n$.
\item $\secat_G(ev) \leq n$: there exist $G$-invariant open sets $U_1,\ldots, U_k$ which cover $X \times X$ and $G$-equivariant sections $s_i :  U_i \rightarrow X^I$ such that $ev \circ s_i$ is  $G$-homotopic to $U_ i \rightarrow X \times X$.

\item $_{\Delta(X)}\cat_G(X \times X) \leq n$: there exist $G$-invariant open sets $U_1,\ldots, U_k$ which cover $X \times X$ which are $G$-compressible into $\Delta(X)$.

\end{enumerate}
\end{thm}

We have the following cohomological lower bound for $\TCG(X)$, using equivariant cohomology theory. Denote by $H^*_G(X)$ the Borel $G$-equivariant cohomology of $X$, with coefficients in an arbitrary commutative ring.

\begin{prop}
Let $\mathfrak{Z}$ be the kernel of the homomorphism $H^*_G(X\times X)\to H^*_G(X)$ induced by the diagonal, then
$\TC_{G} (X) >\nil{(\mathfrak{Z})}$.
\end{prop}


\begin{exam} \label{example1}
\begin{enumerate}
\item  For the free $S^1$-action on itself by rotations, $\TC_{S^1}(S^1)=2$ while $\TC(S^1 / S^1) =\TC(*)=1$, see example 5.10 in \cite{Colman2012}.
\item For the antipodal  action of $\mathbb{Z}_2$   on $S^n$, we have
\[
\TC_{\mathbb{Z}_2}(S^n)= \left\{ \begin{array}{rl} 3 & \mbox{if } n \mbox{ is even}\\
                                        2   & \mbox{if } n \mbox{ is odd}\end{array}\right.
\]  while $\TC(S^n / \mathbb{Z}_2) = \TC(\mathbb{RP}^n) \geq n+1$. See lemma 4.1 in \cite{productprojective} where  equivariant vector fields are used to construct the $G$-equivariant motion planners. In general  $\TC(\mathbb{RP}^n)$ for $n \neq 1,3,7$  is equal to the smallest $k$ such that there is an immersion of $\mathbb{RP}^{n}$ in $\mathbb{R}^{k-1}$, see \cite{Farber2003}.
\item For the reflection  action of $\mathbb{Z}_2$   on $S^n$ with $n\neq 0$,
$$
{\TC}_{\mathbb{Z}_2}(S^n ) \leq  \cat_{\mathbb{Z}_2 \times \mathbb{Z}_2}(S^n \times S^n ) \leq 2 \cat_{\mathbb{Z}_2}(S^n)-1 = 2\times 2 -1=3.
$$
\begin{itemize} 
\item For $n \neq 1$ odd,  ${\TC}_{\mathbb{Z}_2}(S^n)=3$, because $(S^n)^{\mathbb{Z}_2} = S^{n-1}$ and    $3={\TC}(S^{n-1}) \leq {\TC}_{\mathbb{Z}_2}(S^n)$  by proposition \ref{prop3}.
\item For $n$ even,  ${\TC}_{\mathbb{Z}_2}(S^n)=3$, because $3={\TC}(S^{n})$ and ${\TC}(S^{n}) \leq  {\TC}_{\mathbb{Z}_2}(S^n ) $ by proposition \ref{prop3}.
\end{itemize}
We have therefore shown that
\[
\TC_{\mathbb{Z}_2}(S^n) = \left\{ \begin{array}{rl} \infty & \mbox{if } n=1\\
                                        3    & \mbox{if } n\ge 2. \end{array}\right.
\] while $\TC(S^n / \mathbb{Z}_2) = \TC (D^n) =1$.
\item Let $X=S^1 \setminus \{N,S\}$ where $N=(0,1)$ and $S=(0,-1)$.
For the $(\mathbb{Z}_2 \times \mathbb{Z}_2)$-action on $X$ given by the reflection across the $x$-axis and the $y$-axis, we have $\TC_{\mathbb{Z}_2 \times \mathbb{Z}_2}(X) = \infty$. The action has no global fixed points and the space $X$ is $\left ( \mathbb{Z}_2 \times \mathbb{Z}_2 \right )$-compressible into an orbit that is disconnected, see example 2.10 in \cite{BlaszczykKaluba}.
\end{enumerate}
\end{exam}
From the previous examples we see that if $G$ acts freely on $X$, not necessarily $\TC_G (X)={\TC}(X/G)$ and there are no inequalities of the type $\TC_G(X) \leq \TC(X/G)$ or $\TC(X/G) \leq \TC_G(X)$ for general actions.

Examples show that $\TCG(X)$ can be equal to $\TC(X)$, or at the other extreme, one can be finite and the other infinite as shows the case in which  $X$ is a $G$-manifold which is connected but not $G$-connected. 

For a group acting on itself by left translations, we have that $\TCG(G)=\cat(G)$, so that category of topological groups is obtained as a special case of equivariant topological complexity. 

\begin{thm} Let $G$ be a connected metrizable group acting on itself by left translation. Then $\TCG(G)=\cat(G)$. \label{tcgroups}
\end{thm}

Note that this shows that even for spaces that are $G$-compressible to an orbit the equivariant topological complexity can be arbitrarily high.

\begin{prop}
Let $X$ be a $G$-connected topological group. Assume that $G$ acts on $X$ by topological group homomorphisms. Then $\TCG(X)=\cat_G(X)$.
\end{prop}

The next result  relates equivariant and nonequivariant topological complexity for $G$-bundles.
\begin{thm}
Let $E \rightarrow B$ be a numerable principal $G$-bundle and $X$ a $G$-space, then
$$
{\TC}(X_{G}) \leq {\TC}_G(X){\TC}(B)
$$
where $X_G = E \times_G X$ is the total space of  the associated bundle over $B$.
\end{thm}

\begin{thm}\cite{Lubawski2014}
\begin{enumerate}

\item Let $X$ be a $H$-space and $Y$ a $K$-space. Consider $X \times Y$ as a $(H \times K)$-space. Then
$$\TC_{H \times K}(X \times Y) \leq \TC_H(X) + \TC_K(Y)-1.$$
\item From the previous, for $G$-spaces $X$ and $Y$ and $X \times Y$ equipped with the diagonal action, we have that
$$\TC_{G} (X \times Y)\leq \TC_G(X) + \TC_G(Y)-1.$$
See also theorem 4.1 in \cite{productprojective} for the special case of $G$-manifolds.
\end{enumerate}
\end{thm}


\subsection{Invariant topological complexity (Lubawski-Marzantowicz)}

Consider the space $X^I \times_{X/G} X^I= \big\{(\alpha,\beta) \in X^I \times X^I : G\alpha(1)=G\beta(0) \big\} $. The  map $\pi \colon X^I \times_{X/G} X^I\to X \times X$ given by $ \pi(\alpha,\beta) = \big(\alpha(0), \beta(1)\big)$ is a $(G\times G)$-fibration with respect to the obvious actions.

\begin{defn}
The {\em invariant topological complexity} of $X$, $\TC^G(X)$, is the least integer $k$ such that $X \times X$ may be covered by $k$ $(G \times G)$-invariant open sets
  $\{U_1,\ldots ,U_k\}$, on each of which there is a $(G \times G)$-equivariant section 
   $ s_i:{U}_i \to X^I \times_{X/G} X^I$ 
  such that  the  diagram commutes:
$$\xymatrix{& X^I \times_{X/G} X^I \ar[d]^{\pi} \\ {U}_i  \ar[ur]^{s_i} \ar[r]^{}& X\times X }$$
\end{defn}

\begin{prop} Let $X$ and $Y$ be $G$-spaces.
\begin{enumerate}
\item If $X^G \neq \emptyset$, then $X$ is $G$-compressible to a global fixed point if and only if ${\TC}^G(X)=1$. See corollary 2.8 in \cite{BlaszczykKaluba}.
\item ${\TC}^G(X)$ can be finite even when  $X$ is not $G$-connected.
\end{enumerate}
\end{prop}

By restricting a $(G \times G)$-equivariant deformation $c : U \times I \rightarrow X \times X$  to the $(G \times G)$-fixed point set, we have 

\begin{prop}
$\TC(X^G) \leq \TC^G(X).$ \label{prop14}
\end{prop} 

Note that  in general there are no inequalities  of the type  ${\TC}(X^H) \leq \TC^G(X)$ for other subgroups $H \leq G$  as item (5) of example \ref{example} below shows.

Next we show the relationship of invariant topological complexity with equivariant Lusternik-Schnirelmann category.
\begin{prop} Let $x\in X$ and $Gx$ be the orbit of $x$. If $X$ is a $G$-connected space then
$${\TC}^G(X) \leq _{Gx \times Gx} \cat_{G \times G}(X \times X).$$
\end{prop} 
If $x\in X^G$ and $X$ is $G$-connected  then
$$
_{Gx \times Gx} \cat_{G \times G}(X \times X) = _{x \times x} \cat_{G \times G}(X \times X) = \cat_{G \times G} (X \times X)
$$
and therefore
$${\TC}^G(X) \leq \cat_{G \times G} (X \times X).$$

\begin{prop} \cite{BlaszczykKaluba} If $X^G \neq \emptyset$, then
$\cat_G(X) \leq {\TC}^G(X)$.
\end{prop}

Let $\daleth(X)$ be the saturation of the diagonal $\Delta(X)$ with respect to the $(G\times G)$-action.

\begin{thm} For a $G$-space $X$ the following are equivalent:
\begin{enumerate}
\item ${\TC}^G(X) \leq n$.
\item $\secat_{G \times G}(\pi) \leq n$: there exist $(G \times G)$-invariant open sets $U_1,\ldots, U_k$ which cover $X \times X$ and $(G \times G)$-equivariant sections $s_i :  U_i \rightarrow X^I \times_{X/G} X^I$ such that $ev \circ s_i$ is  $(G \times G)$-homotopic to the inclusion $U_ i \rightarrow X \times X$.
\item $_{  \daleth(X)}\cat_{G\times G}(X \times X) \leq n$: there exist $(G \times G)$-invariant open sets $U_1,\ldots, U_k$ which cover $X \times X$ which are $(G\times G)$-compressible into $\daleth(X)$.
\end{enumerate}
\end{thm}

\begin{exam} \label{example}
\begin{enumerate}
\item ${\TC}^{S^1}(S^1)=1$ for the free $S^1$-action on itself by rotations but $\TC(S^1)=2$.
\item For the  antipodal action of $\mathbb{Z}_2$ on $S^n$, we have
${\TC}^{\mathbb{Z}_2}(S^n ) =  \TC(\mathbb{RP}^n)$. 
\item For the reflection  action of $\mathbb{Z}_2$   on $S^n$  with $n\neq 1$, we have
$$
{\TC}^{\mathbb{Z}_2}(S^n ) \leq  \cat_{\mathbb{Z}_2 \times \mathbb{Z}_2}(S^n \times S^n ) \leq 2 \cat_{\mathbb{Z}_2}(S^n)-1 = 2\times 2 -1=3.
$$
If $n$ is odd, ${\TC}^{\mathbb{Z}_2}(S^n)=3$, since  we have that $(S^n)^{\mathbb{Z}_2} = S^{n-1}$ and $3={\TC}(S^{n-1}) \leq {\TC}^{\mathbb{Z}_2}(S^n )$   
 by proposition \ref{prop14}.

\item ${\TC}^{\mathbb{Z}_2}(S^1)=\infty$ for the reflection action of $\mathbb{Z}_2$  on $S^1$.
\item Let $X=S^1 \setminus \{N,S\}$ where $N=(0,1)$ and $S=(0,-1)$.
The $(\mathbb{Z}_2 \times \mathbb{Z}_2)$-action on $X$ given by the reflection across the $x$-axis and the $y$-axis has no global fixed points. The space $X$ is $G$-compressible to an orbit that is disconnected. Therefore $\TC^{\mathbb{Z}_2 \times \mathbb{Z}_2}(X) = 1$ but for one of the copies of $\mathbb{Z}_2$ inside $ \mathbb{Z}_2 \times \mathbb{Z}_2$ we have $\TC(X^{\mathbb{Z}_2})=\TC(S^0)=\infty$. See example 2.10 in \cite{BlaszczykKaluba}.
\end{enumerate}
\end{exam}

The first two and last examples show that in general there are no inequalities of the form $\TC(X) \leq {\TC}^G(X)$ or ${\TC}^G(X) \leq \TC(X)$. Also the last example  shows that there are no inequalites of the form $\TC(X^H) \leq \TC^G(X)$ for $H \leq G$.

The case of ${\TC}^{\mathbb{Z}_2}(S^n)$ for $n$ even is still open. This case is treated  in example 4.2 in \cite{Lubawski2014} but there is a mistake that is explained in remark 3.7 in \cite{BlaszczykKaluba}, it does not follow like in item (3) of example \ref{example1} because in general we do not have $\TC(X) \leq {\TC}^G(X)$.

Since the quotient space of $X \times X$ by the $(G \times G)$-action is $X/G \times X/G$,  we have that $(G \times G)$-invariant open sets in $X \times X$ give open sets in $X/G \times X/G$. Moreover
$(G \times G)$-equivariant local sections of $X^I \times_{X/G} X^I$ give local sections of $(X/G)^I$ and therefore,

\begin{prop}
$\TC(X/G) \leq \TC^G(X).$
\end{prop} 
There is no analogue of this inequality for equivariant topological complexity as the examples in the previous section  show.

For a free action, we can lift a deformation $ c: U \times I \rightarrow X/G \times X/G$ to the diagonal $\Delta(X/G)$, to an equivariant one, by using the Palais covering homotopy theorem. Since  $\pi^{-1}(\Delta(X/G))=\daleth(X)$ under the projection $\pi : X \times X \rightarrow X/G \times X/G$, we have that the lift is a $(G \times G)$-equivariant deformation to $\daleth(X)$. This proves $\TC^G(X) \leq \TC(X/G)$ and gives the following fundamental result.

\begin{prop}
If $G$ is a compact Lie group that acts freely on a metrizable space  $X$ then 
$${\TC}^G (X)={\TC}(X/G).$$
More generally, the result holds if  $X$ has one orbit type. 
\end{prop}

In general it is not known if for a general compact Lie group $\daleth(X) \subset X \times X$ is a $(G \times G)$-cofibration, it is only known for $G$ finite group and $X$ a compact $G$-ANR \cite{Lubawski2014}. This is relevant because of the following product formula.


\begin{prop} Let $X$ be an $H$-space and $Y$ be a $K$-space.
If the inclusion $ \daleth(X) \subset X \times X$ is  a $(H\times H)$-cofibration and $  \daleth(Y) \subset Y \times Y$ is a $(K\times K)$-cofibration, then 
$$
{\TC}^{H \times K}(X \times Y) \leq {\TC}^H(X) + {\TC}^K(Y)-1.
$$
 \end{prop}

 From the previous, for $G$-spaces $X$ and $Y$ we do not obtain the inequality
$${\TC}^{G}(X \times Y) \leq {\TC}^G(X) + {\TC}^G(Y)-1$$
as shown by the example of the free  $S^1$-action on $S^1$ satisfying  ${\TC}^{S^1}(S^1)=1$ and ${\TC}^{S^1}(S^1 \times S^1)=2$. See remark 3.20 in \cite{Lubawski2014}.


\subsection{Strongly equivariant topological complexity (Dranishnikov)}

Let $G\times G$ act on $X\times X$  where $G$ acts in each component.  \textit{Strongly equivariant topological complexity} of $X$, ${\TC}_G^*(X)$, is defined like the equivariant topological complexity, only that  $X\times X$ is now viewed as a $(G\times G)$-space and the open cover is required to be $(G\times G)$-invariant.

\begin{defn}
 The {\em strongly equivariant topological complexity} of $X$,  ${\TC}^*_G(X)$,  is the least integer $k$ such that $X \times X$ may be covered by $k$ $(G \times G)$-invariant open sets
  $U_1,\ldots ,U_k$, on each of which there is a $G$-equivariant section  $ s_i:{U}_i{\to} X^I $ 
for the diagonal action on ${U}_i$, such that  the  diagram commutes:
$$\xymatrix{& X^I  \ar[d]^{\ev} \\ {U}_i  \ar[ur]^{s_i} \ar[r]^{}& X\times X }$$
\end{defn}

\begin{prop} Let $X$ and $Y$ be $G$-spaces.
\begin{enumerate} 
\item For a $G$-connected space $X$  with $X^G \neq \emptyset$, we have that $\TC^*_G(X)=1$ iff $X$ is $G$-contractible.
\item $\TC^*_G(X)=\infty$ if $X$ is not $G$-connected.
\end{enumerate}
\end{prop}

It is obvious that ${\TC}_G(X) \leq  {\TC}_G^*(X)$ 
and therefore using the properties of equivariant topological complexity  we have the following.

\begin{prop} If $H$ is a closed subgroup of $G$, then
\begin{enumerate}
\item ${\TC}(X^H) \leq \TC_G^*(X)$.
\item $\TC_H(X) \leq \TC_G^*(X)$, in particular 
${\TC}(X) \leq \TC_G^*(X)$.
\end{enumerate}
\end{prop}

\begin{prop} For a $G$-connected space with $X^G \neq \emptyset$ we have that
 $$
 {\TC}_G^*(X) \leq \cat_{G \times G} ( X \times X).
 $$ 
\end{prop}

Since ${\TC}_G(X) \leq  {\TC}_G^*(X)$ and $\cat_G(X) \leq \TC_G(X)$  for a $G$-connected space $X$ with $X^G \neq \emptyset$, we have the following,

\begin{prop} For a $G$-connected space with $X^G \neq \emptyset$
$$
\cat_G(X) \leq \TC_G^*(X).
$$
\end{prop}

We can characterize strongly equivariant topological complexity  in terms of relative sectional category and relative Clapp-Puppe category. 

\begin{thm} For a $G$-space $X$ the following are equivalent:
\begin{enumerate}
\item ${\TC}_G^*(X) \leq n$.
\item There exist $(G \times G)$-invariant open sets $U_1,\ldots, U_k$ which cover $X \times X$ and $G$-equivariant sections $s_i :  U_i \rightarrow X^I$ such that $ev \circ s_i$ is  $G$-homotopic to $U_ i \rightarrow X \times X$.
\item There exist $(G \times G)$-invariant open sets $U_1,\ldots, U_k$ which cover $X \times X$  which are $G$-compressible into $\Delta(X)$. 
\end{enumerate}
\end{thm}

\begin{exam}
\begin{enumerate}

\item For a  Lie group acting on itself by translation, we have that ${\TC}_{G}^*(G) \geq \TC_{G}(G)=\cat(G)$ by proposition \ref{tcgroups} but $\TC^G(G)=1$.
\item For the reflection  action of $\mathbb{Z}_2$   on $S^n$, similar arguments as in  example \ref{example1}, prove that 
\[
\TC_{\mathbb{Z}_2}^*(S^n) = \left\{ \begin{array}{rl} \infty & \mbox{if } n=1\\
                                        3    & \mbox{if } n\ge 2. \end{array}\right.
\]

\end{enumerate}
\end{exam}

In case that the group is discrete and the action is free we can describe $\daleth(X)$ as $(G \times G)$-space. By freeness
$$ 
\daleth(X) \cong G \times \Delta(X)
$$
and the action is given by $(h,k) \cdot (g,(x,x)) = (hgk^{-1},(kx,kx))$.

Now restricting to the diagonal subgroup, we have a $G$-action given by 
$$
(hgh^{-1},(hx,hx)).
$$

The quotient $\daleth(X)/G$ is a disjoint union of copies of the diagonal of $X/G$ indexed by the conjugacy classes of $G$, this is a subspace of $X \times_G X$. Since the action is free, there is a covering space $X \times_G X  \rightarrow X/G \times X/G$.

If $X$ is simply connected, the fundamental group of $X \times_G X$ is $G$. Given a deformation $H: U \times I \rightarrow X/G \times X/G$ from the inclusion of an open set to a map $c: U \rightarrow X/G \times X/G$ with $c(U) \subseteq \Delta(X/G)$, the lifting property of covering spaces   gives maps $ \overline{H} : U \times I \rightarrow X \times_G X$ covering $H$. 

Since $\daleth(X)/G$ is a disjoint union of copies of the diagonal of $X/G$ , the lift $\overline{H}(x,t)$ can be chosen so $\overline{H}(x,1)\subseteq \Delta(X/G) \subseteq X \times_G X$. By  using the homotopy lifting property of covering spaces, there is a $G$-equivariant deformation $\overline{H}: \overline{U} \times I \rightarrow X \times X$ that covers $H$ with the property that $\overline{H}(x,1) \subseteq \Delta(X)$. The open set $\overline{U}$ is  $\pi^{-1} (U)$ where $\pi : X \times X \rightarrow X/G \times X/G$ and therefore is $(G \times G)$-invariant. 
\begin{prop} If a  discrete group $G$ acts freely on the simply connected space $X$ then
$$
{\TC}_G^*(X) \leq \TC(X/G).
$$
\end{prop}

The examples above show that this is not true in general.

\begin{thm}
Let $p:X\to B$ be a $F$-bundle between locally compact metric ANR-spaces with the structure group $G$ acting properly on $F$. 
Then $$\TC(X)\le \TC(B)+\TC_{G}^*(F)-1.$$
\end{thm}

\begin{prop}
Suppose that a discrete group $G$ acts freely and properly on a simply connected locally compact ANR-space $Y$. Then $\TC^*_{G}(Y)\le \dim Y+1$. 
\end{prop}

By using the two previous results,

\begin{thm}
Let $X$ be a  CW-complex with fundamental group $G$.  Then
$$
\TC(X)\le \TC(G)+\dim X.
$$
\end{thm}


\subsection{Effective topological complexity (B{\l}aszczyk-Kaluba)}

\begin{defn}
    The {\em effective  topological  complexity of $X$},  $\TC^{G,\infty}(X)$,  is the minimum of the numbers $\TC^{G,n}$ defined as the  least integer $k$ such that
  $X \times X$ may be covered by $k$  open sets
  $U_1,\ldots ,U_k$, on each of which there is a (nonequivariant) section
   $s_i: {U}_i{\to} X^I \times_{X/G} X^I \times_{X/G} \cdots \times_{X/G} X^I =: \mathcal{P}_n(X)$ 
 such that  the  diagram commutes:
$$\xymatrix{& \mathcal{P}_n(X) \ar[d]^{\pi_n} \\ {U}_i  \ar[ur]^{s_i} \ar[r]^{}& X\times X }$$
\end{defn}


\begin{prop} Let $X$ and $Y$ be  $G$-spaces. 
\begin{enumerate}
\item $\TC^{G,n+1}(X) \leq \TC^{G,n}(X)$.
\item $\TC^{G,n}(X)=1$ if $X$ is contractible or $G$-contractible, but the converse is not true.
\item $\TC^{G,n}(X)=\infty$ if $X/G$ is not path connected.
\end{enumerate}
\end{prop}


We have the following cohomological lower bound using non-equivariant cohomology theory of the quotient space.

\begin{prop} Let $\mathfrak{Z}$ be the kernel of zero divisors in the cohomology of $X/G$ with rational coefficients. If $G$ is finite then
$\TC^{G,n} (X) > \nil{( \mathfrak{Z})}$ for $n \geq 2$.
\end{prop}

When the action on the cohomology is trivial, the previous result gives a very useful method to calculate $\TC^{G,n}(X)$.

\begin{prop}
Let $G$ be a finite group acting on $X$ such that the $G$-action on the rational cohomology is trivial. If
$
\TC(X) =\nil{( \mathfrak{Z})}+1
$,
then for $n \geq 2$
$
\TC(X) = \TC^{G,n}(X).
$
\end{prop}
These cohomological bounds do not work with arbitrary coefficients as the last example below shows.


\begin{exam}
\begin{enumerate}
\item  Let $\mathbb{Z}_p$ act on a sphere $S^n$ with $p$ a prime number, $p> 2$.  By using the previous result we have,
$$
\TC^{\mathbb{Z}_p,\infty}(S^n) = 
\left\{ \begin{array}{rl} 3 & \mbox{if } n \mbox{ is even}\\
                                        2   & \mbox{if } n \mbox{ is odd}. \end{array}\right.
$$
\item  If $\mathbb{Z}_2$ acts on a sphere $S^n$ preserving the orientation, the previous result holds.
\item  If $\mathbb{Z}_2$ acts freely on $S^n$, then $\TC^{\mathbb{Z}_2,\infty}(S^n)=2$, while $\TC(S^n / \mathbb{Z}_2) = \TC(\mathbb{RP}^n) \geq n+1$. 
\item For the reflection  action of $\mathbb{Z}_2$   on $S^n$,  $\TC^{\mathbb{Z}_2,\infty}(S^n)=1$.
\item In \cite{BlaszczykKaluba} there are examples of $\mathbb{Z}_p$-actions ($p >2$ prime) on spheres $S^n$ ($n \geq 5)$ with fixed point set $(S^n)^{\mathbb{Z}_p}$ homology spheres that are essential manifolds of dimension $n-2$. The category of a $(n-2)$-dimensional closed essential manifold is $n-1$ by \cite{Berstein} and therefore $\cat_{\mathbb{Z}_p}(S^n) \geq \cat \left( \left (S^n\right )^{\mathbb{Z}_p}\right ) \geq n-1$ while $\TC^{\mathbb{Z}_p,\infty}(S^n) \leq 3$ as in the first example.
\item For a finite group $G$, consider the total space $EG$ of the universal $G$-bundle with its $G$-action. Since the space is contractible, we have that $\TC^{G,\infty}(EG) \leq \TC(EG)=1$, but $\cat_G(EG) = \cat(BG)=\infty$.
\end{enumerate}
\end{exam}

From the previous examples we see that if $G$ acts freely on $X$, not necessarily $\TC^{G,\infty} (X)={\TC}(X/G)$ and that the equivariant Lusternik-Schnirelmann category cannot be a lower bound. See \cite{effective} for a nice discussion.

Since $\TC^{G,1}(X)=\TC(X)$ and $\TC^{G,2}(X) \leq \TC^G(X)$, we have that 
$$\TC^{G,\infty}(X) \leq \min \{\TC(X),\TC^G(X) \}.$$ Together with the fact that for free actions
$\TC^G(X)=\TC(X/G)$, we have

\begin{prop} For a free action of a group $G$ on $X$
$$
\TC^{G,\infty}(X) \leq \TC(X/G).
$$
\end{prop}


\begin{prop} If $H$ is a subgroup of $G$, then
$\TC^{G,n}(X) \leq \TC^{H,n}(X)$.
\end{prop}
\begin{prop} 
Let $X$ be a $H$-space and $Y$ a $K$-space. Then
$$\TC^{H \times K,n}(X \times Y) \leq \TC^{H,n}(X) + \TC^{K,n}(Y)-1.$$
\end{prop}
From the previous, for $G$-spaces $X$ and $Y$ we do not obtain the inequality$$\TC^{G,n}(X \times Y) \leq \TC^{G,n}(X) + \TC^{G,n}(Y)-1$$
as shown by the example of $\mathbb{Z}_2$ acting on $S^1$ by a reflection. From the examples above, we know that $\TC^{\mathbb{Z}_2,\infty}(S^1)=1$ but under the diagonal action $S^1 \times S^1 / \mathbb{Z}_2$ is homeomorphic to $S^2$ and therefore by the cohomological bound stated above $\TC^{\mathbb{Z}_2,\infty}(S^1 \times S^1) > 1$, see section 6 in \cite{effective}.


\section{Comments}

We believe that all the invariants presented in this survey contribute to the study of classical topological complexity and also provide, in different degrees, a valid invariant to further investigate motion planning problems with symmetries.
Much of ordinary topological complexity theory can be adapted and extended to the setting of spaces with an action of a group and each of these notions takes a different approach to that end. 

If we view a $G$-space as being described by the diagram of its fixed  points $X^H$ for the various subgroups $H$ of $G$, the {\em equivariant} and {\em invariant} topological complexities provide a way to import this point of view to the motion planning problems with symmetries. 

From another point of view, we believe that the topological complexity of a $G$-space should reflect in a certain way the {\em quotient object} defined by the action. In particular, in case that the action is free, the {\em invariant} topological complexity of a $G$-space $X$ coincides exactly with the topological complexity of the orbit space $X/G$. 

As a future direction, we envision a notion of topological complexity that recreates a version of this property even in the case when the action is not free. The most suitable object to substitute the orbit space in this case is what is called a  topological {\em orbifold}. 

There has been much recent interest in the study of orbifolds, which 
can locally be described as the quotient of an open subset of Euclidean space by the action of a finite group. In this setting, different group actions may define the same orbifold. Specifically, representable orbifolds are given by a {\em Morita equivalence class} of group actions. If the action is free, it is Morita equivalent to the trivial action on its quotient space. 

We think that a notion of topological complexity for a  $G$-space should display the property that the invariant topological complexity exhibits for free actions but in greater generality. Namely, if two actions are {\em Morita equivalent}, they should have the same topological complexity.  Since the  representation of an orbifold by a group action is not unique, we need a topological complexity that is a true invariant of the orbifold structure and not of the particular representation, i.e. it needs to be checked that we get the same topological complexity for all Morita equivalent actions. We do not know if the invariant topological complexity satisfies this property, but we believe it is the notion best equipped to be developed to be the topological complexity of the actual object defined by the orbits of the action.

\bibliographystyle{amsalpha}

\end{document}